\documentclass[12pt,reqno]{amsart}

\newcommand{\ver}{{\it }}

\usepackage{latexsym}

\newcommand{\szego}{Szeg\"o }

\newcommand{\inv}{^{-1}}
\newcommand{\kahler}{K\"ahler }
\newcommand{\sqrtn}{\sqrt{N}}
\newcommand{\wt}{\widetilde}
\newcommand{\wh}{\widehat}

\newcommand{\R}{{\mathbb R}}
\newcommand{\C}{{\mathbb C}}

\renewcommand{\d}{\partial}
\renewcommand{\le}{\leq}
\newcommand{\dbar}{\bar\partial}

\newcommand{\D}{{\mathbf D}}
\renewcommand{\H}{{\mathbf H}}
\newcommand{\half}{{\frac{1}{2}}}
\newcommand{\vol}{{\operatorname{Vol}}}

\renewcommand{\phi}{\varphi}

\newcommand{\ccal}{\mathcal{C}}

\newcommand{\hcal}{\mathcal{H}}
\newcommand{\lcal}{\mathcal{L}}

\newcommand{\jcal}{\mathcal{J}}

\newcommand{\al}{\alpha}
\newcommand{\be}{\beta}
\newcommand{\ga}{\gamma}
\newcommand{\Ga}{\Gamma}

\newcommand{\la}{\lambda}
\newcommand{\ep}{\varepsilon}
\newcommand{\de}{\delta}
\newcommand{\De}{\Delta}
\newcommand{\om}{\omega}

\newtheorem{theo}{{\sc Theorem}}[section]
\newtheorem{cor}[theo]{{\sc Corollary}}
\newtheorem{lem}[theo]{{\sc Lemma}}

\newenvironment{rem}{\medskip\noindent{\it Remark:\/} }{\medskip}
\newenvironment{defin}{\medskip\noindent{\it Definition:\/} }{\medskip}

\title[Asymptotics of almost holomorphic sections: addendum \ver]
{Asymptotics of almost  holomorphic sections of ample line bundles on
symplectic manifolds: an addendum}

\author{Bernard Shiffman}
\author{Steve Zelditch \ver}
\address{Department of Mathematics, Johns Hopkins University, Baltimore,
MD
21218, USA}
\email{shiffman@math.jhu.edu, zelditch@math.jhu.edu}

\thanks{Research partially supported by NSF grants \#DMS-9800479,
\#DMS-0100474 (first author) and \#DMS-0071358 (second author).}

\date{August 3, 2001}

\begin{document}

\begin{abstract}We define a Gaussian measure on
the space $H^0_J(M, L^N)$ of almost holomorphic sections of powers
of an ample line bundle $L$ over a symplectic manifold $(M,
\omega)$, and calculate the joint probability densities of
sections taking prescribed values and covariant derivatives at a
finite number of points. We prove that they have a universal
scaling limit as $N \to \infty$.  This result is used in
\cite{BSZ2} to extend our previous work \cite{BSZ1} on
universality of scaling limits of correlations between zeros to
the almost-holomorphic setting.

\end{abstract}

\maketitle

\section{Introduction}

This note is an addendum to our study in \cite{SZ2} of almost
holomorphic sections of powers of ample line bundles $L^N \to M$
over almost complex symplectic manifolds $(M, \omega,J)$.
Motivated by the  important role now played by `asymptotically
holomorphic' sections in symplectic geometry (see \cite{A, DON.1}
and many related articles), we studied in \cite{SZ2} a
conceptually related but different  class $H^0_J(M, L^N)$
of `almost holomorphic' sections defined  by a method of Boutet de
Monvel and Guillemin  \cite{BG}. The main
results of \cite{SZ2} were the scaling limit law of the almost-complex \szego
projectors
$$\Pi_N : L^2(M, L^N) \to H^0_J(M, L^N)$$
and various of its geometric consequences.

Our purpose in this addendum is to develop  \cite{SZ2} in a
probabilistic direction,  in the spirit of our earlier work with P. Bleher on
the holomorphic case \cite{BSZ1}, and to complete the discussion in
\cite{BSZ2} of correlations between zeros in the symplectic case.  The space
$H^0_J(M, L^N)$ is finite dimensional and carries a natural Hermitian inner
product (see (\ref{RR}) and (\ref{inner2})). We
endow the unit sphere $SH^0_J(M, L^N)$  in $H^0_J(M,L^N)$ with
Haar probability measure $\nu_N$ and consider the joint probability
distribution (JPD)

$$\D^N_{(z^1,\dots,z^n)}=D^N(x^1,\dots,x^n,\xi^1,\ldots,\xi^n;z^1,\dots,z^n)dxd\xi$$
of a random section
$s_N\in SH^0(M,L^N)$ having
the prescribed values $x^1,\dots,x^n$ and derivatives $\xi^1,\ldots,\xi^n$
at $n$ points
$z^1,\dots, z^n\in M$.  To be more precise, we choose local complex
coordinates $\{z_1,\dots,z_m\}$ and a nonvanishing local
section $e_L$ of $L$ on an open set containing the points $\{z^1,\dots,z^n\}$;
then
$\D^N_{(z^1,\dots,z^n)}$ is the JPD  of the $n(2m+1)$ complex random
variables
\begin{equation}\label{randomvar}\begin{array}{r} x^p=\left.\left\langle
e_L^{*N}, s_N\right\rangle\right|_{z^p},\
\xi^p_q=\big\langle e_L^{*N},N^{-\frac{1}{2}}{\nabla}_{\d/\d z_q}
s_N\big\rangle\big|_{z^p},\
\xi^p_{m+q}=\big\langle e_L^{*N},N^{-\frac{1}{2}}{\nabla}_{\d/\d \bar
z_q} s_N\big\rangle\big|_{z^p},\\[10pt] 1\leq p\leq n, 1\leq q\leq
m\;.\end{array}\end{equation}

As the name implies, almost  holomorphic sections behave in the
large $N$ limit much as do holomorphic sections on complex
manifolds. The main result of this note bears this out by showing
that the JPD in the almost complex symplectic case has the same
universal scaling law  as in the holomorphic case, thereby
finishing  the proof of Theorem 1.2 of \cite{BSZ2} on universal
scaling limits of correlations between zeros in the symplectic
case.

\begin{theo}\label{usljpd-sphere} Let $L$ be the complex line bundle over a
$2m$-dimensional compact integral symplectic manifold $(M,\om)$
with curvature $\omega.$  Let $P_0\in M$, and choose a local frame $e$ for
$L$ and local complex coordinates centered at $P_0$ so that $\om|_{P_0}$ and
$g|_{P_0}$ are the usual Euclidean \kahler form and metric
respectively, $\|e_L\|_{P_0} =1$, and $\nabla e_L|_{P_0} = 0$. Then
$$ \D^N_{(z^1/\sqrtn,\dots, z^n/\sqrtn)} \longrightarrow \D^\infty
_{(z^1,\dots,z^n)}\quad\quad
\left((z^1,\dots,z^n)\in\C^{mn}\right)\;,$$   where the measures
$\D^\infty_{(z^1,\dots,z^n)} =
\ga_{\De^\infty(z)}$ are the same universal Gaussian measures
  as in the holomorphic case; in particular,
they are supported on the holomorphic 1-jets. \end{theo}

Let us review the formula for $\ga_{\De^\infty(z)}$. As we will
recall,  a (complex) Gaussian measure on $\C^k$ is a measure of the
form
\begin{equation}\label{cxgauss}\ga_\De= \frac{e^{-\langle\De\inv
z,\bar z\rangle}}{\pi^{k}{\det\De}} dz\,,\end{equation}
where $dz$ denotes Lebesgue measure on $\C^k$, and $\De$ is a
positive definite Hermitian $k\times k$ matrix.  The  matrix
$\De=\big(\De_{\al\be}\big)$ is the covariance matrix of $\ga_\De$:
\begin{equation}\label{covar}
\langle
z_\al\bar z_\be\rangle_{\ga_\De} =\De_{\al\be},\quad\quad 1\leq \al,\be\leq
k\;.\end{equation}
Since the universal limit measures $\ga_{\De^\infty(z)}$ vanish along
non-holomorphic directions, they are singular measures on the
space of all 1-jets. To deal with singular measures, we introduce in
\S
\ref{s-gaussians} generalized Gaussian measures whose covariance matrices
(\ref{covar}) are only semi-positive definite; a generalized Gaussian is
simply a Gaussian measure supported on the subspace corresponding to the
positive eigenvalues of the covariance matrix.

The covariance matrix ${\De^\infty(z)}$ is given along holomorphic
directions by the same formula as in the holomorphic case \cite[(97)]{BSZ1},
namely
\begin{equation}\label{usldelta} \De^\infty(z)= \frac{m!}{c_1(L)^m}\left(
\begin{array}{cc}
A^\infty(z) & B^\infty(z) \\
B^{\infty}(z)^* & C^\infty(z)
\end{array}\right)\;,\end{equation} where $z=(z^1,\dots,z^n)$ and
\begin{eqnarray*}
 A^\infty(z)^p_{p'} &=& \Pi_1^\H(z^p,0;z^{p'},0) \,,\quad\quad
\Pi_1^\H(u,0;v,0) =
\frac{1}{\pi^m} e^{u \cdot\bar{v} - \half(|u|^2 + |v|^2)}\;,\\
B^\infty(z)^{p}_{p'q'} &=&\left\{\begin{array}{ll}
(z^p_{q'}-z^{p'}_{q'}) \Pi_1^\H(z^p,0;z^{p'},0) \
&\mbox{for}\quad  1\le
q\le m\\ 0 & \mbox{for}\quad  m+1\le q\le 2m\end{array}\right.\ ,\\
C^\infty(z)^{pq}_{p'q'} &=&\left\{\begin{array}{ll}
\left(\delta_{qq'}+(\bar z^{p'}_q -\bar z^p_q)
(z^p_{q'}-z^{p'}_{q'})\right)\Pi_1^\H(z^p,0;z^{p'},0) \
&\mbox{for}\quad  1\le q,q'\le m\\ 0 & \mbox{for}\quad
\max(q,q')\ge m+1\end{array}\right.\ .
\end{eqnarray*}

In other words, the coefficients of $\De^\infty(z)$ corresponding to the
anti-holomorphic directions vanish, while the coefficients corresponding to
the holomorphic directions are given by
the \szego kernel $\Pi_1^\H$ for the reduced Heisenberg group (see
\cite[\S 1.3.2]{BSZ1}) and its covariant derivatives.

A technically interesting novelty in the proof is
the role of the $\bar{\partial}$ operator. In the holomorphic case,
$\D^N_{(z^1,\dots,z^n)}$ is supported
on the subspace of jets of sections satisfying $\dbar s=0$.  In the
almost
complex case,  sections do not satisfy this equation, so
$\D^N_{(z^1,\dots,z^n)}$ is a measure on the  higher-dimensional space of
all 1-jets. However, Theorem \ref{usljpd-sphere} says that the mass in the
non-holomorphic directions shrinks to zero as $N \to \infty$.

An alternate statement of Theorem \ref{usljpd-sphere} involves equipping
$H^0_J(M, L^N)$ with a Gaussian measure, and letting
$\wt\D^N_{(z^1,\dots,z^n)}$ be the corresponding joint probability
distribution, which is a Gaussian measure on the complex
vector space of 1-jets of sections.  We show (Theorem \ref{usljpd}) that
these Gaussian measures
$\wt\D^N$ also have the same scaling limit $\D^\infty$, so that
asymptotically the probabilities are the same as in the holomorphic case,
where universality was established in
\cite{BSZ2}.   It is then easy to see that
$\wt\D^N_{(z^1,\dots,z^n)}=\ga_{\De^N}\,$ where $\De^N$ is the covariance
matrix of the random variables in (\ref{randomvar}).   The main step in
the proof is to show that the covariance matrices $\De^N$ underlying $\wt\D^N$
tend in the scaling limit to a semi-positive matrix $\De^{\infty}$.  It
follows that the scaled distributions $\wt\D^N$ tend to the generalized
Gaussian
$\ga_{\De^{\infty}}$ `vanishing in the $\bar{\partial}$-directions.'

 In a subsequent article \cite{SZ3}, we  obtain further probabilistic results
on holomorphic and almost holomorphic sections. Regarding almost holomorphic
sections,  we prove that a sequence $\{s_N\}$ of almost holomorphic sections
satifies the bounds
$$\|s_N\|_{\infty}/\|s_N\|_{L^2} = O(\sqrt{\log N}), \;\;\;\;
\|\bar{\partial} s_N\|_{\infty}/\|s_N\|_{L^2} = O(\sqrt{\log
N})$$ almost surely. Hence almost holomorphic
sections satisfy similar bounds to asymptotically holomorphic sections in
\cite{DON.1, A}.

Finally, we mention some intriguing questions relating our probabalistic approach to almost holomorphic
sections to the now-standard complexity-theoretic approach to asymptotically holomorphic sections in symplectic
geometry, due to Donaldson and further developed by
Auroux and others. From an analytical viewpoint (which of course is just one technical side of their work), the
key results are existence
theorems for one or several asymptotically holomorphic sections satisfying
quantitative transversality conditions, such as
$$ s(z) = 0 \ \Longrightarrow\|\bar{\partial} s(z) | < |\partial
s(z)|\quad\quad
\forall z \in M$$ in the case of one section. Can one use the probabalistic
method to prove such existence results?   It is the global nature of the
problem which makes it difficult. On small balls, our methods rather easily
give  lower bounds for quantitative transversality of the type:
$$\mu\{s: \ |\bar{\partial} s(z) | < |\partial s(z)|\;\;\; \forall z \in
B_{\frac{D}{\sqrt{N}}}(z_0) \; s. th. \; s(z) = 0\ \} \ >\  1-
\frac{C_{\ep}}{N^{1-\ep}}\,.$$ However,  there are roughly $C_m
N^m$ balls of radius $1/\sqrt{N}$, so one cannot simply sum this
small-ball estimate.  To globalize, one would need to partition
$M$ into small cubes as in \cite{DON.1} and then analyze the
dependence of transversality conditions from one cube to another.

\begin{rem} This paper and \cite{SZ2} were  originally contained  in the
archived  preprint \cite{SZ1}.  In revising that paper for publication as
\cite{SZ2}, we expanded the section on the \szego
kernel, added a discussion on transversality, and relocated the
material on the JPD to this article. In particular,  Theorem
\ref{usljpd-sphere}  is Theorem 0.2 of \cite{SZ1}, and Theorem \ref{usljpd} is
Theorem 5.4 of
\cite{SZ1}.  \end{rem}

\section{Background}\label{background}

To avoid duplication, we only set up some basic notation and background and
refer to \cite{BSZ1, BSZ2, SZ2} for further discussion and details.

Let  $(M, \omega, J)$ denote  a compact almost-complex symplectic manifold such that
$[\frac{1}{\pi}\omega]$ is an integral
cohomology class, and  such that
$\om(Jv,Jw)=\om(v,w)$ and
$\omega(v, Jv) > 0$.
We further let
 $(L, h, \nabla) \to M$ denote a quantizing Hermitian line bundle and   metric connection with
 curvature $\frac{i}{2}\Theta_L = \omega$.
We denote by $L^N$ the
$N^{\rm th}$ tensor power of $L$.

 We give $M$  the Riemannian metric
$g(v,w)=\om(v,Jw)$. We   denote by
$T^{1,0}M,
$ resp.\
$T^{0, 1}M$,  the holomorphic, resp.\ anti-holomorphic, sub-bundle of the
complex tangent bundle $TM$;  i.e., $J = i$ on $T^{1,0}M$ and $J = -i$ on
$T^{0,1}M$.

We now recall our notion of `preferred coordinates' from \cite{SZ2}.  They are
important because the universal scaling laws are only valid in
such coordinates.  A coordinate system $(z_1,\dots,z_m)$ on a
neighborhood $U$ of a point $P_0\in M$ is {\it preferred\/} at $P_0$  if
any two of the following conditions (and
hence all three) are satisfied:

\begin{enumerate}
\item[i)] $\quad\d/\d z_j|_{P_0}\in T^{1,0}(M)$, for $1\leq j\leq m$,
\item[ii)] $\quad\om({P_0})=\om_0$,
\item[iii)] $\quad g({P_0} )= g_0$,
\end{enumerate}
where $\omega_0$ is the standard
symplectic form and $g_0$ is the Euclidean metric:
$$\om_0=\frac{i}{2}\sum_{j=1}^m dz_j\wedge d\bar z_j =\sum_{j=1}^m
(dx_j\otimes dy_j - dy_j\otimes dx_j)\,,\quad g_0=
\sum_{j=1}^m
(dx_j\otimes dx_j + dy_j\otimes dy_j)\,.$$
Here we write
$z_j=x_j+iy_j$, and we let
$\{\frac{\d}{\d z_j},
\frac{\d}{\d\bar z_j}\}$ denote the dual frame to $\{dz_j, d\bar z_j\}$.

As  in \cite{BSZ1, BSZ2, SZ2}, it is advantageous to work on the associated
principal
$S^1$ bundle $X \to M$, and our \szego kernels will be defined there.  Let
$\pi: L^* \to M$ denote the dual line bundle to $L$ with dual metric $h^*$, and
put $X = \{v \in L^*:
\|v\|_{h^*} =1\}$.   We let $r_{\theta}x =e^{i\theta} x$ ($x\in X$) denote the
$S^1$ action on $X$.  We then identify  sections $s_N$ of $L^N$ with
equivariant function
$\hat{s}$ on $L^*$ by the rule
\begin{equation} \label{sNhat}\hat{s}_N(\lambda) = \left( \lambda^{\otimes N},
s_N(z)
\right)\,,\quad
\la\in X_z\,,\end{equation} where $\lambda^{\otimes N} = \lambda \otimes
\cdots\otimes
\lambda$; then $\hat s_N(r_\theta x) = e^{iN\theta} \hat s_N(x)$. We denote by
$\lcal^2_N(X)$ the space of such equivariant functions transforming by the
$N^{\rm th}$ character.

When working on $X$, covariant derivatives on sections of $L$ become
horizontal derivatives of equivariant functions.
We consider preferred coordinates  $(z_1,\dots,z_m)$ centered at a point
$P_0\in M$ and  a local frame $e_L$ for $L$ such that $\|e_L\|_{P_0}
=1$ and
$\nabla e_L|_{P_0} = 0$.  This gives us
coordinates $(z_1,\dots,z_m,\theta)$ on $X$  corresponding to
$x=e^{i\theta}\|e_L(z)\| e_L^*(z)\in X$. We showed in \cite[\S
1.2]{SZ2} that
\begin{equation}\frac{\d^h}{\d z_j} = \frac{\d}{\d z_j}
+\left[\frac{i}{2}\bar z_j +O(|z|^2)\right]
\frac{\d}{\d \theta}\,,\quad \frac{\d^h}{\d\bar z_j} = \frac{\d}{\d\bar z_j}
-\left[\frac{i}{2} z_j +O(|z|^2)\right]
\frac{\d}{\d \theta}
\,,\label{dhoriz}\end{equation}
 where $\frac{\d^h}{\d
z_j}$ (resp.\ $\frac{\d^h}{\d
\bar z_j}$) denotes the horizontal lift of
$\frac{\d}{\d z_j}$ (resp.\ $\frac{\d}{\d
\bar z_j}$).

The almost-complex Cauchy-Riemann
operator $\dbar_b$ on $X$ does not satisfy $\bar{\partial}_b^2 = 0$ in
general and usually has no kernel.
Following a method of  Boutet de Monvel -Guillemin, we defined
in \cite{SZ2} the space ${\mathcal H}_J^N$ of equivariant almost-CR
functions on $X$ as the kernel of a certain deformation $\bar{D}_0$ of the
$\bar{\partial}_b$ operator on $\lcal^2_N(X)$. The space $H^0_J(M, L^N)$ is
then the corresponding space of sections. The
\szego kernel $\Pi_N(x, y)$ is the kernel of the orthogonal projection
$\Pi_N: \lcal^2_N(X) \to {\mathcal H}_J^N.$ The dimension
$d_N = \dim H_J^0(M, L^N)$ satisfies the Riemann-Roch formula (see \cite{BG})
\begin{equation} \label{RR} d_N  = \frac{c_1(L)^m}{m!}N^m + O(N^{m-1})
\;.\end{equation}
Since $H_J^0(M, L^N)$ is finite dimensional, the \szego
kernel $\Pi_N$ is smooth and is given by:
$$\Pi_N(x, y) = \sum_{j=1}^{d_N} S_j^N(x)\overline{ S_j^N(y)}\;,$$
where $\{S_j^N\}$ is an orthonormal basis for $H^0_J(M,L^N)$.

It would take
us too far afield to discuss the definition or significance of the spaces
$H^0_J(M, L^N)$ here; we refer the reader to \cite{SZ2} for background.

\section {A generalized Poincar\'e-Borel lemma}\label{s-gaussians}

In this section, we give a generalization of the Poincar\'e-Borel lemma
concerning the asymptotics of linear push forwards of measures on spheres of
growing dimensions, which we shall use in the proof of Theorem
\ref{usljpd-sphere}.

Recall that a Gaussian measure on $\R^n$ is a measure of the form
$$\ga_\De = \frac{e^{-\half\langle\De\inv
x,x\rangle}}{(2\pi)^{n/2}\sqrt{\det\De}}  dx_1\cdots dx_n\,,$$ where $\De$
is a positive definite symmetric $n\times n$ matrix.  The matrix $\De$ gives
the second moments of $\ga_\De$:
\begin{equation}\label{moments}\langle x_jx_k
\rangle_{\ga_\De}=\De_{jk}\,.\end{equation} This
Gaussian measure is
also characterized by its Fourier transform
\begin{equation}\label{muhat}\wh{\ga_\De}(t_1,\dots,t_n) = e^{-\half\sum
\De_{jk}t_jt_k}\,.\end{equation}
If we let $\De$ be the $n\times n$ identity matrix, we obtain the standard
Gaussian measure on $\R^n$,
$$\ga_n:=\frac{1}{(2\pi)^{n/2}}e^{-\half |x|^2} dx_1\cdots dx_n\,,$$
with the
property that the
$x_j$ are independent Gaussian variables with mean 0 and variance 1.

By a {\it generalized complex Gaussian measure\/} on $\C^n$, we
mean a generalized Gaussian measure $\ga_\De^c$ on $\C^n\equiv
R^{2n}$ with moments
$$\big\langle z_j\big\rangle_{\ga_\De^c} =0,\quad \big\langle z_j z_k
\big\rangle_{\ga_\De^c}=0,\quad \big\langle z_j\bar
z_k\big\rangle_{\ga_\De^c} =\De_{jk},\quad\quad 1\leq j,k\leq n,$$
where $\De=(\De_{jk})$ is an $n\times n$ positive semi-definite
Hermitian matrix; i.e. $\ga_\De^c=\ga_{\half \De^c}$, where
$\De^c$ is the $2n\times 2n$ real symmetric matrix of the inner
product on $R^{2n}$ induced by $\De$. As we are interested here in
complex Gaussians, we  drop the `$c$' and write
$\ga_\De^c=\ga_\De$.  In particular, if $\De$ is  a strictly
positive Hermitian matrix, then $\ga_\De$ is given by
(\ref{cxgauss}).

The push-forward of a Gaussian measure by a surjective linear map is also
Gaussian. In the next section, we shall push forward Gaussian
measures on the spaces
$H^0_J(M,L^N)$ by linear maps that are sometimes not surjective.  Since
these non-surjective push-forwards are  singular measures, we need to
consider the case where
$\De$ is positive semi-definite.  In this case, we use (\ref{muhat}) to
define a measure $\ga_\De$, which we call a {\it generalized Gaussian\/}.
If $\De$ has null eigenvalues, then the generalized Gaussian
$\ga_\De$ is a Gaussian measure on the subspace $\Lambda_+\subset\R^n$
spanned by the positive eigenvectors.  (Precisely, $\ga_\De=\iota_*
\ga_{\De|\Lambda_+}$, where $\iota:\Lambda_+\hookrightarrow \R^n$ is the
inclusion.  For the completely degenerate case
$\De=0$, we have
$\ga_\De=\de_0$.)  Of course, (\ref{moments}) also holds for semi-positive
$\De$. One useful property of generalized Gaussians is that the push-forward
by a (not necessarily surjective) linear map $T:\R^n\to\R^m$ of a
generalized Gaussian $\ga_\De$ on $\R^n$
is a generalized Gaussian on $\R^m$:
\begin{equation}\label{pushgaussian} T_*\ga_\De=\ga_{T\De T^*}\end{equation}
Another useful property of  generalized Gaussians is the following fact:

\begin{lem}\label{continuity}  The map $\De\mapsto\ga_\De$ is a continuous
map from the positive semi-definite matrices to the space of positive
measures on
$\R^n$ (with the weak topology).\end{lem}

\begin{proof} Suppose that $\De^N\to\De^0$. We must show that
$(\De^N,\phi)\to(\De^0,\phi)$ for a compactly supported test function
$\phi$.  We can assume that $\phi$ is $\ccal^\infty$. It then
follows from (\ref{muhat}) that
$$(\ga_{\De^N},\phi)=(\wh{\ga_{\De^N}},\wh\phi)\to
(\wh{\ga_{\De^0}},\wh\phi)=(\ga_{\De^0},\phi)\,.$$\end{proof}

We shall use the following `generalized Poincar\'e-Borel lemma' relating
spherical measures to Gaussian measures in our proof of Theorem
\ref{usljpd-sphere} on asymptotics of the joint probability distributions for
$SH^0_J(M,L^N)$.

\begin{lem}\label{spherical-vs-gaussian}
 Let $T_N: \R^{d_N} \to R^k$, $N=1,2,\dots$, be a sequence of linear maps,
where $d_N
\to \infty$.  Suppose that $\frac{1}{d_N}T_N T_N^*
\to
\Delta$. Then
$T_{N*} \nu_{d_N} \to \gamma_{\Delta}$. \end{lem}

\begin{proof}
Let $V_N$ be a $k$-dimensional subspace of $\R^{d_N}$ such that
$V_N^\perp\subset \ker T^N$, and let
$p_N:\R^{d_N}\to V_N$ denote the orthogonal projection.  We decompose
$T_N=B_N\circ A_N$, where
$A_N=d_N^{1/2}p_N:\R^{d_N}\to V_N$, and $B_N =
d_N^{-1/2}T_N|_{V_N}:V_N\buildrel{\approx }\over\to\R^k$.
  Write $$A_{N*}\nu_{d_N} =
\al_N\,,\quad\quad T_{N*}
\nu_{d_N} =
B_{N*}\al_N = \be_N\,.$$  We easily see that (abbreviating $d=d_N$)

\begin{equation}\label{archimedes}\al_N=A_{N*} \nu_d = \psi_d dx\,,\quad\quad
\psi_d =
\left\{\begin{array}{cl}\frac{\sigma_{d-k}}{\sigma_d d^k} [1 -  \frac{1}{d}
|x|^2]^{(d-k-2)/2}
\ &{\rm for}\ |x|<\sqrt{d}\\[6pt]
\ 0 & {\rm otherwise}\end{array}\right.,\end{equation}
where $dx$ denotes Lebesgue measure on $V_N$, and
$\sigma_n=\vol(S^{n-1})=
\frac{2\pi^{n/2}}{\Ga(n/2)}$. (The case $k=1, \
d=3$ of (\ref{archimedes}) is  Archimedes' formula \cite{Arc}.) Since
$[1 -   |x|^2/d]^{(d-k-2)/2}
\to e^{-|x|^2/2}$ uniformly on compacta and $\frac{\sigma_{d-k}}{\sigma_d
d^k}\to \frac{1}{(2\pi)^{k/2}}$, we conclude that $\al_N\to\ga_k$. (This is
the Poincar\'e-Borel Theorem; see Corollary \ref{PB} below.) Furthermore,
$$\left(1 -  \frac{1}{d} |x|^2\right)^{(d-k-2)/2}
\leq \exp\left(-\frac{d-k-2}{2d}|x|^2\right)\leq
e^{\frac{k+2}{2}}e^{-\half|x|^2}\quad\quad{\rm
for}\ d\ge k+2\,,\
|x|\leq \sqrt{d}\,,$$
 and hence \begin{equation}\label{ldc}\psi_{d_N}(x)\leq
C_k e^{-|x|^2/2}\,.\end{equation}

Now let $\phi$ be a compactly supported continuous test function on
$\R^k$. We must show that \begin{equation}\label{true}\int \phi
d\be_N\to\int
\phi d\ga_{\De}\,.\end{equation}  Suppose on the contrary that
(\ref{true}) does not hold. After passing to a subsequence, we may assume that
$\int \phi d\be_N\to c\ne \int
\phi d\ga_{\De}$.  Since the eigenvalues of $B_N$ are bounded, we can
assume  (again taking a subsequence) that $B_N\to B_0$, where
$$B_0 B_0^*=\lim _{N\to\infty}B_NB_N^*=\lim
_{N\to\infty}\frac{1}{d_N} T_N T_N^*=\De\,.$$  Hence,
$$\int_{\R^k} \phi
d\be_N =\int_{V_N} \phi(B_Nx)\psi_{d_N}(x)dx \to \int_{V_N} \phi(B_0
x)\frac{e^{-|x|^2/2}}{(2\pi)^{k/2}} dx=
 \int_{V_N} \phi(B_0 x)d\ga_k(x)\,,$$
where the limit holds by dominated convergence, using (\ref{ldc}).
By (\ref{pushgaussian}), we have $B_{0*}\ga_k=\ga_{B_0B_0^*}=\ga_\De$, and
hence
$$ \int_{V_N} \phi(B_0 x)d\ga_k(x) = \int_{\R^k}\phi d\ga_{\De}\,.$$
Thus (\ref{true}) holds for the subsequence, giving a contradiction.
\end{proof}

\medskip We note that the above proof
began by establishing the {\it Poincar\'e-Borel Theorem\/} (which is
a special case of the  of Lemma \ref{spherical-vs-gaussian}):

\begin{cor} {\rm  (Poincar\'e-Borel)} Let $P_d:\R^d\to \R^k$ be given by
$P_d(x)=\sqrt{d}(x_1,\dots,x_k)$.  Then $$P_{d*}\nu_d \to \ga_k\,.$$
\label{PB}\end{cor}

\section {Proof of Theorem \ref{usljpd-sphere}}{\label{s-jpd}
We return
to
our complex Hermitian line bundle $(L,h)$ on a compact almost complex
$2m$-dimensional symplectic manifold
$M$
with symplectic form $\omega=\frac{i}{2}\Theta_L$, where $\Theta_L$ is the
curvature of $L$ with respect to a connection $\nabla$.
We now describe the $n$-point joint distribution arising from our
probability space $(SH^0_J(M,L^N),\nu_N)$.
We introduce  the
Hermitian inner product on $H^0_J(M,L^N)$:
\begin{equation}\label{inner2}\langle s_1, s_2 \rangle = \int_M h^N(s_1,
s_2)\frac{1}{m!}\om^m \quad\quad (s_1, s_2 \in
H^0_J(M,L^N)\,)\;,\end{equation}
and we write
$\|s\|_2=\langle s,s \rangle^{1/2}$. Recall that $SH^0_J(M,L^N)$ denotes the
unit sphere $\{\|s\|=1\}$ in $H^0_J(M,L^N)$ and $\nu_N$ denotes its Haar
probability measure.

We let $J^1(M,L^N)$ denote the space of
1-jets of sections of
$L^N$. Recall that we have the exact sequence of vector bundles
\begin{equation}\label{jet} 0\to T^*_M\otimes L^N   \stackrel{\iota}{\to}
J^1(M,L^N) \stackrel{\rho}{\to}
L^N\to 0\,.\end{equation} We consider the jet
maps
$$J^1_z:H^0_J(M,L^N)\to J^1(M,V)_z\,,\quad J^1_zs=\ \mbox{the
1-jet
of}\ s\ \mbox{at}\ z\,, \quad \mbox{for}\ z\in M\,.$$
The covariant derivative $\nabla:J^1(M,L^N)\to T^*_M\otimes
L^N$  provides a splitting
of (\ref{jet}) and an  isomorphism
\begin{equation}\label{splitting}(\rho,\nabla):J^1(M,L^N){\buildrel{\approx}
\over
\longrightarrow}  L^N\oplus(T^*_M\otimes L^N)\,.\end{equation}

\begin{defin} \label{DEFJPD}
The {\it $n$-point joint probability distribution\/} at
points $P^1,\dots, P^n$ of $M$ is the probability measure
\begin{equation}\label{alternately}
\D^N_{(P^1,\dots,P^n)}:=(J^1_{P^1}\oplus\cdots
\oplus J^1_{P^n})_*\nu_N\end{equation} on the space
$J^1(M,L^N)_{P^1}\oplus\cdots \oplus J^1(M,L^N)_{P^n}$.
\end{defin}

Since we are interested in the scaling limit of $\D^N$, we need to describe
this measure more explicitly:  Suppose that
$P^1,\dots,P^n$ lie in a coordinate neighborhood of a point $P_0\in M$ and
choose  preferred coordinates $(z_1,\dots,z_m)$  at $P_0$. We let
$z^p_1,\dots,z^p_m$ denote the coordinates of the point $P^p$ ($1\leq p\leq
n$), and we write
$z^p=(z^p_1,\dots,z^p_m)$. (The coordinates of
$P_0$ are $0$.) We consider the
$n(2m+1)$  complex-valued random
variables $x^p,\ \xi^p_q$ ($1\leq p\leq n,\ 1\leq q\leq 2m$)
on $S\hcal^2_N(X)\equiv SH^0_J(M,L^N)$ given by
\begin{equation}\label{dms1def} x^p(s) = s(z^p,0)\,,
\end{equation}  \begin{equation}
\label{dms2def} \xi^p_q(s)=\frac{1}{\sqrtn} \frac{\d^h s}{\d
z_q}(z^p)\,,\quad
\xi^p_{m+q}(s)=\frac{1}{\sqrtn} \frac{\d^h s}{\d\bar z_q}(z^p)\quad\quad (1\leq
q\leq m)\,,
\end{equation}
for $s\in SH^0_J(M,L^N)$.

We now write $$x=(x^1,\dots,x^p)\,,\quad
\xi=(\xi^p_q)_{1\leq p\leq n,1\leq q\leq 2m}\,,\quad  z=(z^1,\dots,z^n)\,.$$
Using (\ref{splitting}) and the variables $x^p,\ \xi^p_q$ to make the
identification
\begin{equation}\label{identification}
J^1(M,L^N)_{P^1}\oplus\cdots \oplus J^1(M,L^N)_{P^n}\equiv
\C^{n(2m+1)}\,,\end{equation} we can write
$$\D^N_z = D^N(x,\xi;z)dxd\xi\,,$$ where
$dxd\xi$ denotes Lebesgue measure on $\C^{n(2m+1)}$.

Before proving Theorem \ref{usljpd-sphere} on the scaling limit of $\D^N_z$,
we state and prove a corresponding result replacing $(SH^0_J(M,L^N),\nu_N)$
with the essentially equivalent Gaussian space $H^0_J(M,L^N)$ with the
{\it normalized standard Gaussian measure}
\begin{equation}\label{gaussian}\mu_N:= \left(\frac{d_N}{\pi}\right)^{d_N}
e^{-d_N|c|^2}dc\,,\quad\quad s=\sum_{j=1}^{d_N}c_jS_j^N\,,\end{equation}
where
$\{S_j^N\}$ is an orthonormal basis for $H^0_J(M,L^N)$.  This measure  is
characterized by the property that the $2d_N$ real variables
$\Re c_j, \Im c_j$ ($j=1,\dots,d_N$) are independent Gaussian random variables
with mean 0 and variance $1/2d_N$; i.e.,
\begin{equation}\label{normalized}\langle c_j \rangle_{\mu_N}= 0,\quad
\langle c_j c_k\rangle_{\mu_N} = 0,\quad \langle c_j \bar c_k
\rangle_{\mu_N}=
\frac{1}{d_N}\de_{jk}\,.\end{equation} Our normalization guarantees that the
variance of
$\|s\|_2$ is unity: $$\langle
\|s\|^2_2\rangle_{\mu_N}=1\,.$$ We then consider the {\it Gaussian joint
probability distribution\/}
\begin{equation}
\wt \D^N_{(P^1,\dots,P^n)}= \wt D^N(x,\xi;z)dxd\xi =(J^1_{P^1}\oplus\cdots
\oplus J^1_{P^n})_*\mu_N\,.\end{equation}
Since $\mu_N$ is Gaussian and the map
$J^1_{P^1}\oplus\cdots\oplus J^1_{P^n}$
is linear, it follows  that the joint probability
distribution
is a generalized Gaussian measure of the form
\begin{equation}\label{Dgaussian}D^N(x,\xi;z)dxd\xi=
\ga_{\De^N(z)}\,.\end{equation}
We shall see below that the covariance
matrix $\De^N(z)$ is given in terms of the \szego kernel.

We have the following alternate form of Theorem \ref{usljpd-sphere}:

\begin{theo}\label{usljpd} Under the hypotheses and notation of Theorem
\ref{usljpd-sphere}, we have
$$\wt \D^N_{(z^1/\sqrtn,\dots, z^n/\sqrtn)} \longrightarrow \D^\infty
_{(z^1,\dots,n^n)}\;.$$
\end{theo}
\smallskip\begin{proof}  We use the coordinates $(z_1,\dots,z_m,\theta)$ on
$X$ given by preferred coordinates at $P_0\in M$ and a local frame $e_L$
for $L$ with $\|e_L\|_{P_0}
=1$ and
$\nabla e_L|_{P_0} = 0$ as in \S \ref{background}. The
covariance matrix
$\De^N(z)$ in (\ref{Dgaussian}) is a positive semi-definite $n(2m+1)\times
n(2m+1)$ Hermitian matrix. If the map $J^1_{z^1}\oplus\cdots\oplus J^1_{z^n}$
is surjective, then $\De^N(z)$ is strictly positive definite and $\wt
D^N(x,\xi;z)$ is a smooth function. On the other hand, if the map is not
surjective, then $\wt D^N(x,\xi;z)$ is a distribution supported on a linear
subspace.  For example, in the integrable holomorphic case, $\wt D^N(x,\xi;z)$
is supported on the holomorphic jets, as follows from the discussion below.

By (\ref{moments}), we have
\begin{eqnarray}&
\Delta^N(z)=\left(
\begin{array}{cc}
A & B \\
B^* & C
\end{array}\right)\,,\nonumber \\
&A=\big( A^{p}_{p'}\big)=
\big \langle x^p \bar x^{p'}\big\rangle_{\mu_N}\,,\quad
B=\big(B^{p}_{p'q'}\big)=
\big\langle x^p \bar \xi^{p'}_{q'}\rangle_{\mu_N}\,,\quad
C=\big(C^{pq}_{p'q'}\big)=
\big\langle  \xi^{p}_{q}\bar \xi^{p'}_{q'}\rangle_{\mu_N}\,,\label{dms3}\\
& p,p'=1,\dots,n, \quad q,q'=1,\dots,2m\,.\nonumber
\end{eqnarray}
(We note that $A,\ B,\ C$ are $n\times n,\ n\times 2mn,\
2mn\times 2mn$ matrices, respectively; $p,q$ index the rows, and
$p',q'$ index the columns.)

We now describe the the entries of the matrix $\De^N$ in terms of the
\szego kernel.
We have by (\ref{normalized}) and (\ref{dms3}), writing
$s=\sum_{j=1}^{d_N}c_j S^N_j$,
\begin{equation}\label{dms4}A^{p}_{p'}
=\big \langle x^p \bar x^{p'}\big\rangle_{\mu_N}
= \sum_{j,k=1}^{d_N}\big\langle c_j \bar c_k\big\rangle_{\mu_N}
S_j^N(z^p,0)\overline{S_k^N(z^{p'},0)}=
\frac{1}{d_N}\Pi_N(z^p,0;z^{p'},0)\,.\end{equation}
We need some more notation to describe the matrices $B$ and $C$: Write
$$\nabla_q=\frac{1}{\sqrtn}\frac{\d^h}{\d z_q}\,,\quad
\nabla_{m+q}=\frac{1}{\sqrtn}\frac{\d^h}{\d\bar z_q}\,,\quad 1\leq q\leq m\,.$$
We  let $\nabla^1_q$, resp.\ $\nabla^2_q$, denote
the differential operator on $X\times X$ given by applying $\nabla_q$ to the
first, resp.\ second, factor ($1\leq q\leq 2m$).  By differentiating
(\ref{dms4}), we obtain
\begin{eqnarray}\label{dms5} B^{p}_{p'q'} &=&  \frac{1}{d_N}
\overline{\nabla}^2_{q'}\Pi_N(z^p,0;z^{p'},0)\,,\\
\label{dms6}
C^{pq}_{p'q'} &=& \frac{1}{d_N}
\nabla^1_q\overline{\nabla}^2_{q'}\Pi_N(z^p,0;z^{p'},0)\,.\end{eqnarray}

We now use the  scaling asymptotics of the almost holomorphic
\szego kernel $\Pi_N(x,y)$  given in \cite{SZ2}.
In addition to the above assumption on the the local frame $e_L$, we further
assume that $$\quad \nabla^2 e_L|_{P_0} = -(g+i\om)\otimes e_L|_{P_0}\in
T^*_M\otimes T^*_M\otimes L\;.$$
(In
\cite{SZ3}, we called such an $e_L$ a `preferred frame' at $P_0$, and the
resulting coordinates were called  `Heisenberg coordinates.')  Then we have
(see
\cite{SZ3}, Theorem 2.3):
\begin{equation}\label{SCALING} \begin{array}{l}
N^{-m}\Pi_N(\frac{u}{\sqrtn},\frac{\theta}{N};
\frac{v}{\sqrtn},\frac{\phi}{N})\\ \\ \quad = \frac{1}{\pi^m}
e^{i(\theta-\phi)+u \cdot\bar{v} - \half(|u|^2 + |v|^2)}\left[1+ \sum_{r = 1}^{K}
N^{-\frac{r}{2}} b_{r}(P_0,u,v) + N^{-\frac{K +1}{2}}
R_K(P_0,u,v,N)\right],\\ \\ \quad
\mbox{\it where}\ \
\|R_K(P_0,u,v,N)\|_{\ccal^j(\{|u|+|v|\leq \rho\})}\leq C_{K,j,\rho} \ \ \mbox{\it
for}\
\  j=1,2,3,\dots.\end{array}\end{equation}
It follows from  (\ref{dms4})--(\ref{dms6})
and (\ref{SCALING}), recalling (\ref{dhoriz})--(\ref{RR}) that
\begin{equation} \Delta^N(\frac{z}{\sqrtn})\to
\De^\infty(z)= \frac{m!}{c_1(L)^m}\left(
\begin{array}{cc}
A^\infty(z) & B^\infty(z) \\
B^{\infty}(z)^* & C^\infty(z)
\end{array}\right)\end{equation}
in the notation of (\ref{usldelta}).

Finally,  we apply Lemma
\ref{continuity} to (\ref{Dgaussian}) and conclude that
$$\wt\D^N_{z/\sqrtn}=\ga_{\De^N(z/\sqrtn)} \to
\ga_{\De^\infty(z)}=\D^\infty_{z}\,.$$
\end{proof}

\medskip\noindent {\it Proof of Theorem \ref{usljpd-sphere}:\/} The proof is
similar to that of Theorem \ref{usljpd}. This time we define
$$\De^N= \frac{1}{d_N} \jcal_N \jcal_N^*:H^0(M,L^N)\to \C^{n(2m+1)}\,,$$
where $\jcal_N = J^1_{P^1}\oplus\cdots\oplus J^1_{P^n}$ under the
identification (\ref{identification}).  We see immediately that $\De^N$ is
given by (\ref{dms4})--(\ref{dms6}) and the conclusion follows from Lemma
\ref{spherical-vs-gaussian} and (\ref{usldelta}). \qed

\begin{rem} There are other similar ways to define the joint probability
distribution that have the same universal scaling limits.  One of these is to
use the (un-normalized) standard Gaussian measure
$\ga_{2d_N}$ on
$H^0_J(M,L^N)$ in place of the normalized Gaussian $\mu_N$ in Theorem
\ref{usljpd} to obtain joint
densities $D^N_\#(x,\xi;z)=D^N(\frac{x}{N^{m/2}},\frac{\xi}{N^{m/2}};z)$.
Then we would have instead
$$D^N_\#(N^{m/2}x,N^{m/2}\xi;N^{-1/2}z)dxd\xi\to \ga_{\De^\infty(z)}\,.$$
Another similar result is to let $\la_N$ denote normalized Lebesgue measure on
the unit ball $\{\|s\|\leq 1\}$ in $H^0_J(M,L^N)$ and to let
$\wh\D^N_z=\jcal_{N*} \la_N$.  By a similar argument as above, we also have
$\wh\D^N_{z/\sqrtn}\to \ga_{\De^\infty(z)}$.
\end{rem}


\begin{thebibliography}{HHHH}

\bibitem[Arc]{Arc} Archimedes, {\it On the Sphere and Cylinder\/} (Greek),
Syracuse, ca.\ 257BC.



\bibitem[Aur]{A} Denis Auroux,
 Estimated transversality in symplectic geometry and projective maps,  to
appear in Proc. International KIAS Conference (Seoul, 2000),
http://xxx.lanl.gov/abs/math.SG/0010052.

\bibitem[BSZ1]{BSZ1}  P. Bleher, B. Shiffman and S. Zelditch,  Universality
and scaling of correlations between zeros on complex manifolds, {\it
Invent.\ Math.}  142 (2000), 351--395.

\bibitem[BSZ2]{BSZ2}  P. Bleher, B. Shiffman and S. Zelditch, Universality
and scaling of zeros on symplectic manifolds, in {\it Random Matrix Models
and Their Applications}, P. Bleher and A.  Its (Eds.), MSRI Publications
40, Cambridge Univ.\ Press, 2001, http://xxx.lanl.gov/abs/math-ph/0002039.


\bibitem[BSZ3]{BSZ3}  P. Bleher, B. Shiffman and S. Zelditch,  Correlations
between zeros and supersymmetry, {\it Commun.\ Math.\
Phys.}, to appear, http://xxx.lanl.gov/abs/math-ph/0011016.

\bibitem[BoGu]{BG} L.  Boutet de Monvel and V.  Guillemin, {\it The
Spectral Theory of Toeplitz Operators}, {\it Ann.\ Math.\ Studies\/} 99,
Princeton Univ.\ Press, Princeton, 1981.

\bibitem[BoSj]{BS} L. Boutet de Monvel and J. Sj\"ostrand, Sur la
singularit\'e des noyaux de Bergman et de Szeg\"o, {\it Asterisque\/} 34--35
(1976), 123--164.

\bibitem[Don]{DON.1} S. K.  Donaldson, Symplectic submanifolds and almost
complex
geometry, {\it J. Diff.\  Geom.}  44 (1996), 666--705.

\bibitem[ShZe1]{SZ1}  B. Shiffman and S. Zelditch, Random almost holomorphic
sections of ample line bundles on symplectic manifolds, (preprint 2000),
http://xxx.lanl.gov/abs/math.SG/0001102.

\bibitem[ShZe2]{SZ2} B. Shiffman and S. Zelditch, Asymptotics of almost
holomorphic sections of ample line bundles on symplectic manifolds,  {\it J.
Reine Angew.\ Math.}, to appear.

\bibitem[ShZe3]{SZ3} B. Shiffman and S. Zelditch, Random polynomials and Levy
concentration of measure, (in preparation).



\bibitem[Woo]{W} N. M. J. Woodhouse, {\it Geometric Quantization},
Clarendon Press,
Oxford, 1992.

\end{thebibliography}
\end{document}